\documentclass[12pt,a4paper]{article}

\usepackage[english]{babel}
\usepackage{amsmath}
\usepackage{amssymb}
\usepackage{amsthm}
\usepackage{amsfonts}
\usepackage{array}
\usepackage{bm}
\usepackage{booktabs}
\usepackage{chngpage}
\usepackage{epsfig}
\usepackage{float}
\usepackage{graphicx}
\usepackage{lineno}
\usepackage{mdframed}
\usepackage{microtype}
\usepackage{multirow}
\usepackage{subfig}
\usepackage{textcomp}
\usepackage{tabularx}
\usepackage{caption}
\usepackage{threeparttable}
\usepackage{url}
\usepackage{verbatim}
\usepackage{algorithm}
\usepackage{algorithmic}

\graphicspath{{figures/}}

\begin{document}

\title{FT-GCR: a fault-tolerant generalized conjugate residual elliptic solver}
\author{\hspace{-.5cm}Mike Gillard$^{(1)}$, Tommaso Benacchio$^{(2)}$}
\date{}

\maketitle

\begin{center}

{\small
$^{(5)}$ Loughborough University, UK\\
{\tt m.gillard@lboro.ac.uk}
}

\vskip 0.25cm
{\small 
$^{(2)}$ MOX -- Modelling and Scientific Computing, \\
Dipartimento di Matematica, Politecnico di Milano, Milan, Italy \\
{\tt tommaso.benacchio@polimi.it}
}
\end{center}

\vspace*{0.5cm}

\pagebreak

\abstract{With the steady advance of high performance computing systems featuring smaller and smaller hardware components, the systems and algorithms used for numerical simulations increasingly contend with disruptions caused by hardware failures and bit-levels misrepresentations of computing data. In numerical frameworks exploiting massive processing power, the solution of linear systems often represents the most computationally intensive component. Given the large amount of repeated operations involved, iterative solvers are particularly vulnerable to bit-flips. 

A new method named FT-GCR is proposed here that supplies the preconditioned Generalized Conjugate Residual Krylov solver with detection of, and recovery from, soft faults. The algorithm tests on the monotonic decrease of the residual norm and, upon failure, restarts the iteration within the local Krylov space. Numerical experiments on the solution of an elliptic problem arising from a stationary flow over an isolated hill on the sphere show the skill of the method in addressing bit-flips on a range of grid sizes and data loss scenarios, with best returns and detection rates obtained for larger corruption events. The simplicity of the method makes it easily extendable to other solvers and an ideal candidate for algorithmic fault tolerance within integrated model resilience strategies.}

\pagebreak

\section{Introduction}

Numerical models used in the simulation of physical processes rely on dependable hardware as well as on software that is robust against bit-level faults. While both aspects have been taken for granted in the past decades, the emergence and exponential growth of available computing power on the same time frame has spurred a race towards more accurate and efficient simulations, with top high-performance computing facilities fast approaching 1 exaflop and millions of computing nodes \cite{schulthess:2018}. A case in point is numerical weather prediction, where global atmospheric models currently run at resolutions only accessible to local area models a few years ago \cite{dueben:2020}.

As cores multiply in number and shrink in size, distributed computing systems become more vulnerable to loss of resources caused by hardware failures, while algorithms are disrupted by increasingly frequent soft faults, spurious corruptions in their binary representation. Standard disk checkpointing and restart resilience strategies clearly cease to be viable at the point where then overhead associated with the backup procedure exceeds the mean time between failures \cite{giraud:2020}.

At the end point of partially implicit time discretisation strategies for evolution problems, solvers for linear systems and associated preconditioners often represent an efficiency bottleneck in numerical models. The ability of linear solvers to detect, correct and recover from faults with minimal overhead therefore takes on paramount importance in efforts towards integrated model resilience \cite{benacchio:2021}. In this context, Krylov solvers have been equipped with interpolation-restart mechanisms against hard faults without overhead in the case of no failures \cite{agullo:2016}, and multigrid-based compressed checkpointing strategies \cite{goeddeke:2015} have also been proposed. A fault-tolerant version of the GMRES method (FT-GMRES) has been tested in the context of selected reliable and unreliable model sections \cite{hoemmen:2011, sao:2013}. Ultimately, FT-GMRES both relies on an additional flexibility of the underlying adapted model (flexible GMRES) to adopt differing preconditioning and stopping criteria between successive Krylov subspace minimisation steps (outer iterations), and uses a bounded state function (Hessenberg matrix) to detect possible errors (soft-faults). Such options are not directly transferable to other Krylov subspace solvers, where the solution proposed here is generally compatible, but provides a much less robust detection mechanism than that of FT-GMRES.     

This paper proposes a framework for algorithmic fault tolerance and local data recovery within iterative methods used for the solution of linear systems. A Fault-Tolerant Generalised Conjugate Residual (FT-GCR) elliptic solver is designed and implemented with the capability to detect and recover from soft faults. A lightweight mechanism for the detection of soft faults in the form of bit flips within the GCR Krylov solver is introduced, based on a simple check exploiting the monotonic decrease of the solver residual norm at each iteration of the algorithm. In addition, the method is equipped with fault recovery that enables resuming operation following the fault. 

Numerical experiments are carried out on the solution of an elliptic problem arising from the discretisation of steady potential flow past a hill on the sphere. A range of configurations are simulated, with spherical octahedral Gaussian grids up to O640 (corresponding to approximately 15 km horizontal resolution) running on up to 3240 MPI processes, and with various data loss values and injection probabilities of soft fault events, with between one in five to four in a million grid points affected by the soft fault. With the added fault-tolerance, the solver is able to recover effectively from soft faults, reaching convergence in fewer iterations compared with the unprotected solver version. In addition, higher detection rates and return on fault tolerance scores are found for larger corruption events. Finally, the paper contains an evaluation and discussion of the findings with an outlook to future development and applications.

\section{The GCR and a Fault-Tolerant extension}

We consider the following elliptic problem for the unknown $\phi$:
\begin{equation}
\mathcal L\left(\phi\right)=\mathcal R,
\label{eq:BVP}
\end{equation} 
where $\mathcal L$ denotes the elliptic operator and $\mathcal R$ the right-hand side \cite{smolarkiewicz:2000, smolarkiewicz:2005, smolarkiewicz:2011}. 

\noindent The GCR elliptic solver belongs to a class of Krylov subspace methods, minimises the $l_2$ norm of the residual, $r$, of problem \eqref{eq:BVP}, and solves the $k^{th}$-order damped oscillator equation with preconditioner operator $\mathcal P$ \cite{thomas:2003}:
\begin{equation}
\dfrac{\partial^k\mathcal P(\phi)}{\partial\tau^k}+\dfrac{1}{T_{k-1}\left(\tau\right)}\dfrac{\partial^{k-1}\mathcal P(\phi)}{\partial\tau^{k-1}} +...+\dfrac{1}{T_{1}\left(\tau\right)}\dfrac{\partial^1\mathcal P(\phi)}{\partial\tau^1}=\mathcal L\left(\phi\right)-\mathcal R.
\label{eq:damped}
\end{equation}
The damped oscillator equation \eqref{eq:damped} is discretised in pseudo-time $\tau$ and optimal parameters $T_1,..T_{k-1}$ are determined to assure the minimisation of $\left<r,r\right>$ for the field $\phi$, see \cite{smolarkiewicz:2000,smolarkiewicz:2011} for a detailed discussion. The solver employs bespoke left-preconditioning (see, e.g., \cite{kuehnlein:2019}) to address the large condition number induced by the domain anisotropy for global atmospheric problems. 

The iterative nature of the GCR method provides ample opportunity to recompute faulty data at minimal cost to the solver, given an ability to detect such faults. A basic detection method can be derived, given that for each iteration of GCR ($n$ in Algorithm \ref{ftgcrk}), the $k^{th}$ order dampening takes place on a Krylov subspace, $\mathcal K_n$, such that $\mathcal K_{n}\subseteq \mathcal K_{n+1}$, $\forall n\in1,2,\ldots$. As such, each iteration of GCR minimises the $l_2$ norm of the residual, $r$, on that Krylov subspace. Since the norm is non-increasing on each subspace, it is also non-increasing between subspaces. If during computation the discrete $l_2$ norm value increases, that most likely indicates a problem with the solver. 

Algorithm \ref{ftgcrk} describes the fault-tolerant GCR (FT-GCR) method, including the backup and fault detection steps, following \cite{smolarkiewicz:2000} for notation. FT-GCR detection and correction additions are shown in red. The `last known good' solution, denoted $[\#]^\star$, is backed up at after a complete pass of GCR ($n$) where no faults were detected. 

Resiliency testing of GCR($k$) indicates that the method is quite resilient to soft fault events (e.g., bit flips), although such events often cause the solver to stagnate for an iteration. Most often the solver continues on, albeit often converging at a slow rate. This remains the case even for very large fault events (where many data points are corrupted). Fault events sometimes cause convergence issues at much later solver passes, and the response of the solver to fault events is difficult to predict.

\begin{algorithm}[ht]
\footnotesize
\begin{algorithmic}
\STATE For any initial guess, $\phi^0$, set $r^0=\mathcal L\left(\phi^0\right) -\mathcal R$, $p^0=\mathcal P^{-1}\left( r^0\right)$; then iterate: \
\FOR{ $n=1,2,...$  until convergence}  
\FOR{ $\nu=0,...,k-1$} 
\STATE $\beta=\dfrac{\left< r^\nu\mathcal L\left( p^\nu\right)\right>}{\left< \mathcal L\left(p^\nu\right)\mathcal L\left(p^\nu\right)\right>}$  \\
\STATE $\phi^{\nu+1}=\phi^\nu+\beta p^\nu$   \\
\STATE $r^{\nu+1}=r^\nu+\beta\mathcal L\left(p^\nu\right)$ \\
\IF{ $\| r^{\nu+1}\|_2 \le \epsilon$} \STATE exit \ENDIF
\textcolor{red}{\IF{ $ \|r^{\nu+1}\|_2\ge \|r^{\nu}\|_2$} 
\STATE $n=n-1$
\STATE reset $\left[\phi,r,p,\mathcal L(p)\right]^0$ to $\left[\phi,r,p,\mathcal L(p)\right]^\ast$   
\ELSIF { $\nu=0$ }
\STATE set $\left[\phi,r,p,\mathcal L(p)\right]^\ast$ to $\left[\phi,r,p,\mathcal L(p)\right]^0$ \\
\ENDIF} \\
\STATE $e=\mathcal P^{-1}\left( r^{\nu+1} \right)$ \\
\STATE Compute $\mathcal L(e)$ \
\FOR {$l=0,..\nu$}  
\STATE $\alpha_l=\dfrac{\left< \mathcal L(e) \mathcal L\left( p^l\right)\right>}{\left< \mathcal L\left(p^l\right)\mathcal L\left(p^l\right)\right>}$ \\ 
\STATE $p^{\nu+1}=e+{\displaystyle\sum_{l}^{\nu}}\alpha_l p^l$ \  
\STATE $\mathcal L\left(p^{\nu+1}\right)=\mathcal L(e)+{\displaystyle\sum_{l}^{\nu}}\alpha_l\mathcal L\left(p^l\right)$ \\
\ENDFOR \\
\ENDFOR \\
reset $\left[\phi,r,p,\mathcal L(p)\right]^k$ to $\left[\phi,r,p,\mathcal L(p)\right]^0$  
\ENDFOR
\end{algorithmic}
\caption{ FT-GCR($k$) (including red section) - GCR($k$) (excluding red section)}\label{ftgcrk}
\end{algorithm}

Therefore, resilience is built into the GCR solver as follows. When a fault is detected, the GCR solver can be easily reverted to a state previously deemed good (i.e., no fault detected), at the cost of at most a full iteration of Algorithm \ref{ftgcrk} (if a fault occurs during the $\nu=k-1$ pass). The backup of the known good configuration occurs following the exit check and fault test, at the start of a new Krylov subspace loop (i.e., backing up the output of the full $k$ passes on the last Krylov subspace). This ensures the backup solution underwent the full iterative process for a given subspace -- backing up a known good configuration for $\nu<k$ would be akin to varying the order of the damped oscillator \eqref{eq:damped} between outer GCR iterations - which can cause some instability to the overall solution. 

\section{Test problem and numerical results}
\label{sec:num_res}

We consider the case of potential flow $\mathbf{v}$ over a Gaussian-shaped hill, with governing equations:

\begin{equation}\label{eq:pot_flow}
\begin{split}
 & \mathbf{v}=\mathbf{v_a}-\nabla\phi     \\
 & \nabla\cdot\left(\rho\mathbf{v}\right)=0
 \end{split}
\end{equation}
\noindent where $\mathbf{v_a}$ is the ambient velocity and $\rho$ is the prescribed density. The problem is cast in non-orthogonal, terrain-following coordinates, whereby the vertical coordinate is adjusted to the hill profile, depending on the vertical extension of the domain $H$ and bottom topography $h$. The terrain-following reference frame involves metric terms, see  \cite{smolarkiewicz:1994} for details.

Model \eqref{eq:pot_flow} is then discretized in space using a finite volume method in the horizontal and a finite difference method in the vertical. The underlying Octahedral global meshes are generated using the Atlas library \cite{Deconinck:2017}. Combining the transformed equations brings to an elliptic problem of the form \eqref{eq:BVP}, where $\mathcal{L}$ is a discrete Poisson operator, 
which is then solved using the fault-tolerant preconditioned GCR(k) (Algorithm \ref{ftgcrk}). For further details about the system of equations, the GCR solver, and the finite volume discretisation used here, see, e.g., \cite{kuehnlein:2019,smolarkiewicz:2000,smolarkiewicz:2005,smolarkiewicz:2011, smolarkiewicz:1994}. 

The behaviour of FT-GCR in solving the elliptic problem \eqref{eq:BVP} is investigated to determine the response of the fault-tolerant solver at scale. 
MPI-enabled runs are carried out for a range of grid resolutions with suitable core counts and tolerances for the residual $l_2$ norm-based solver exit criterion, and for $0.0004\%$, $0.04\%$, $0.2\%$, $1\%$, $5\%$, and $20\%$ total data loss (Table \ref{tab:resol}). For the simulations in this paper, a purely horizontal ambient velocity with $\mathbf{v_{a,x}}=20\,\textrm{m/s}$ is used, and the values $H=40800\,\textrm{m}$, $h=4000\,\textrm{m}$ and $h_r=3E05\,\textrm{m}$ for domain height and maximum hill height and hill radius are considered.

\begin{table}[ht]
\footnotesize
\centering
\caption{\footnotesize 
Approximate horizontal grid resolutions at the Equator $\Delta x_{EQ}$, tolerance \texttt{tol} on the $l_2$ norm of the residual, number of total grid points, number of MPI processes and grid points per process for the numerical experiments with Octahedral grids O$N$ considered in this paper. A fixed number of 51 vertical levels is used throughout the simulations.}
 \begin{tabular}{lrrrrr}
 \toprule
 \midrule
O$N$              & O40 & O80 & O160 & O320 & O640\\
\midrule
$\Delta x_{EQ}$ [km] & 227 & 119 & 61 & 31 & 15.5 \\
\texttt{tol} & E-14 & E-12 & E-10 & E-8 & E-6 \\  
Points & 399840 & 1452480 & 5516160 & 21477120 & 84733440 \\
MPI procs          & 36 & 108 & 216 & 864 & 3240 \\
Points/proc  & 1110 & 13448 & 25538 & 24858 & 26152 \\  
\bottomrule
\end{tabular}
\label{tab:resol}
\end{table}

Since the simulations are distributed over many computational processes, faults are set to occur on a single MPI process per fault injection event. The faulty process is determined randomly, and the procedure is re-randomised following each injection event. Injection occurs after the preconditioning stage with implausibly high probabilities of either $2\%$ or $5\%$. A maximum of $10$ individual fault events are allowed. The fault causes a given number of entries of $e=P^{-1}\left(r^{\nu+1}\right)$, the preconditioner output, to suffer a bit flip. 
The amount of data entries corrupted during each fault event is varied between runs, with the number of data values corrupted by each event given as a percentage of the total number of array entries. 

\begin{figure}[htbp]
\centering
\subfloat{
\includegraphics[width=0.46 \textwidth]{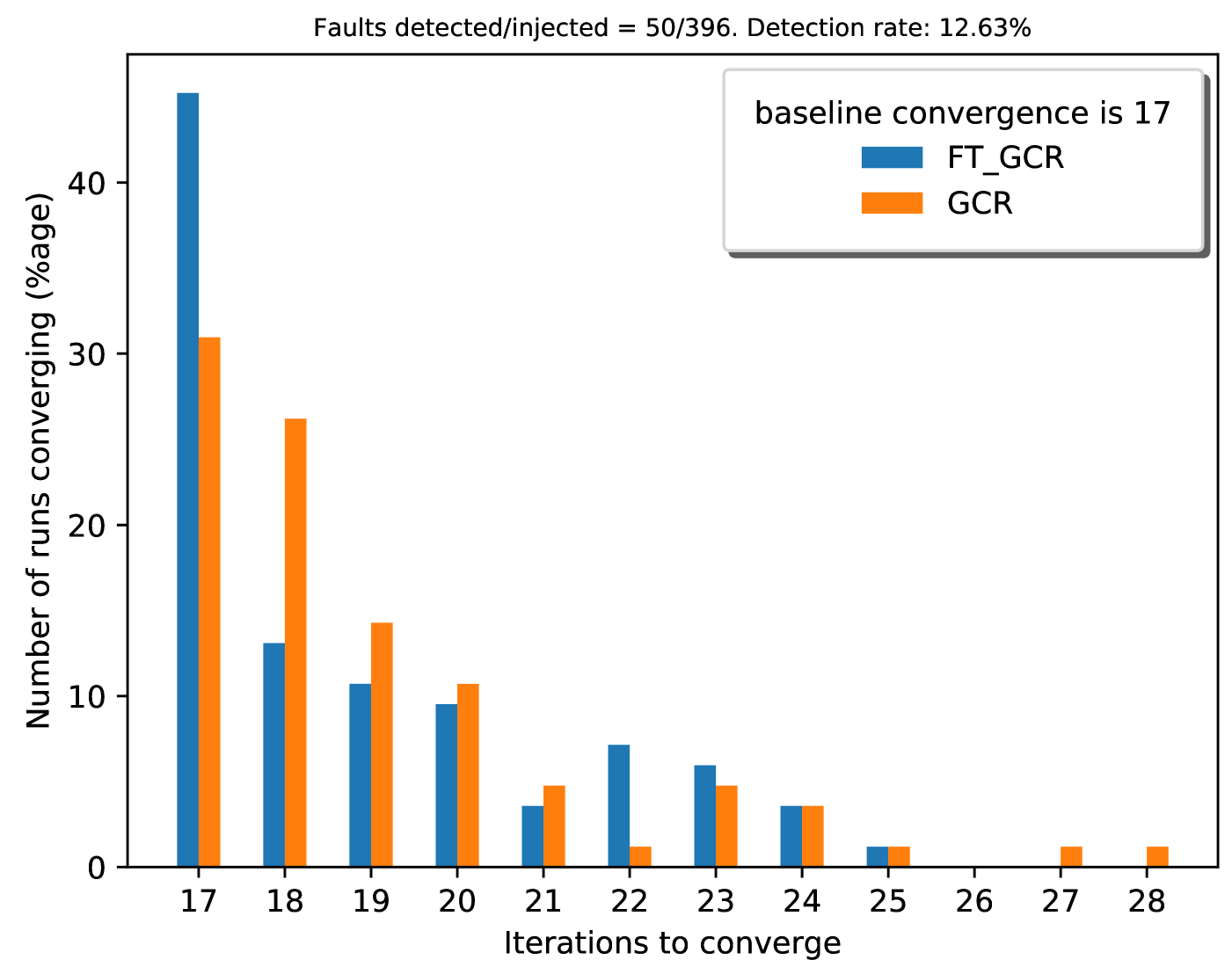}}\subfloat{
\includegraphics[width=0.46 \textwidth]{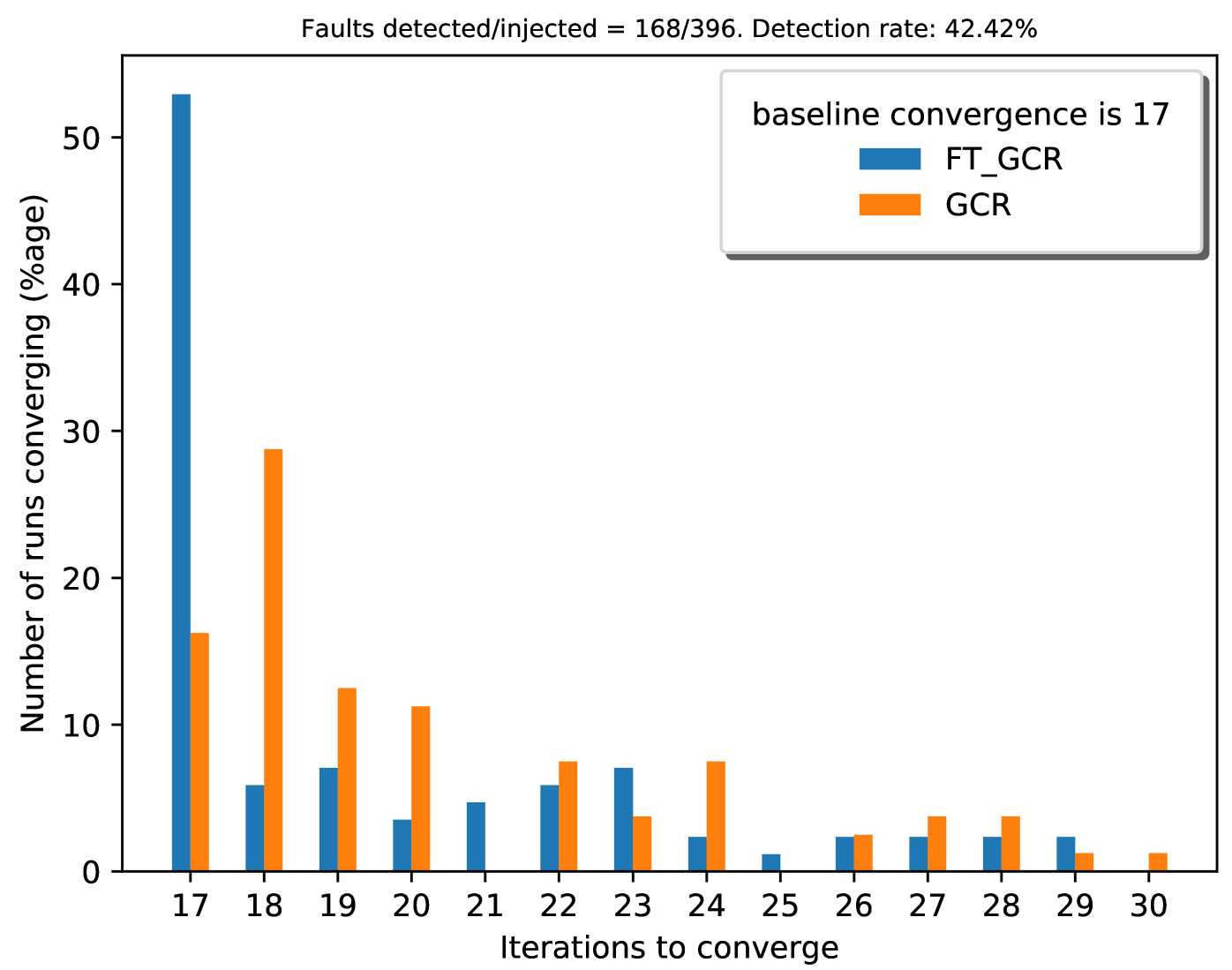}
}\\   
\subfloat{
\includegraphics[width=0.46 \textwidth]{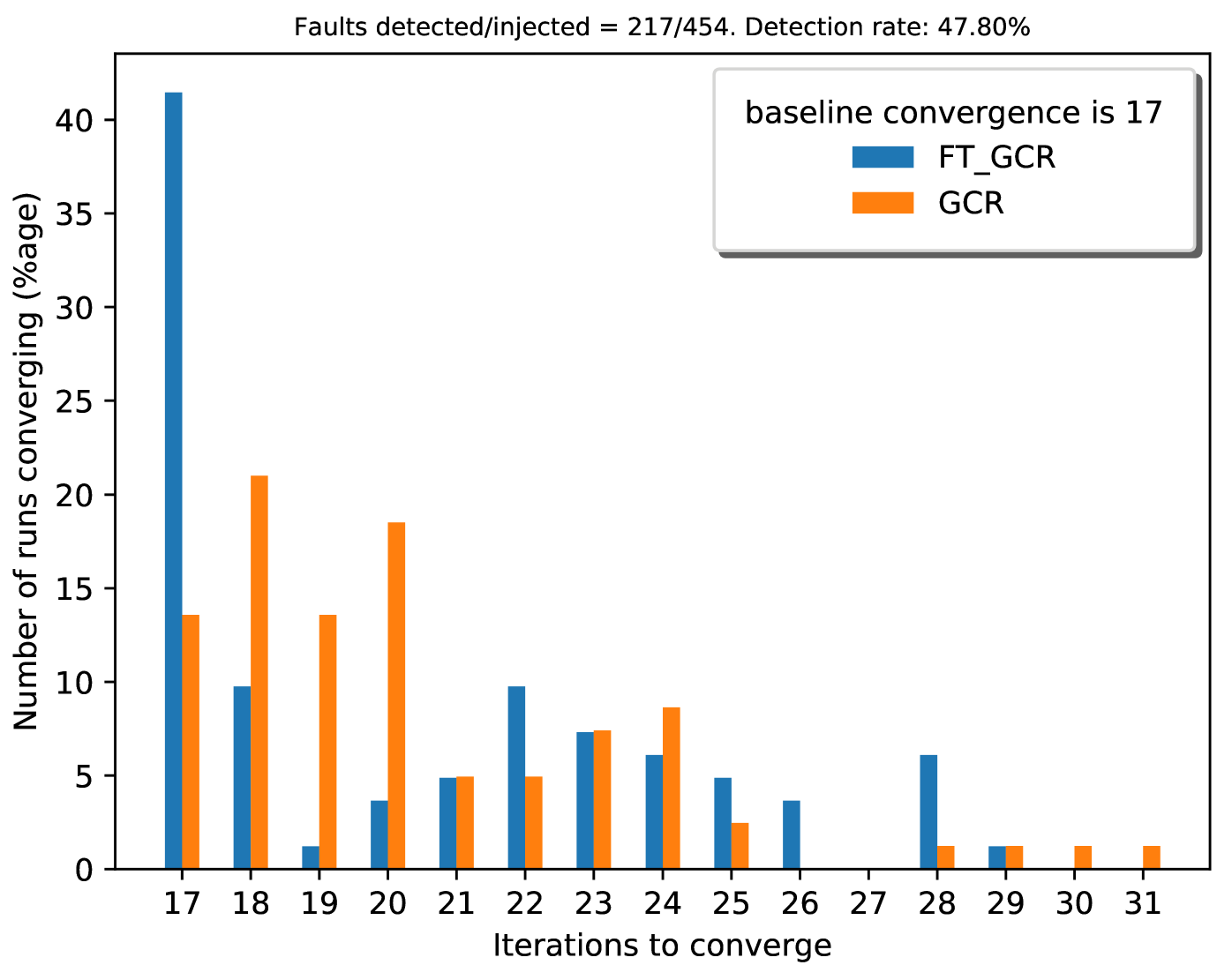}
}\subfloat{
\includegraphics[width=0.46 \textwidth]{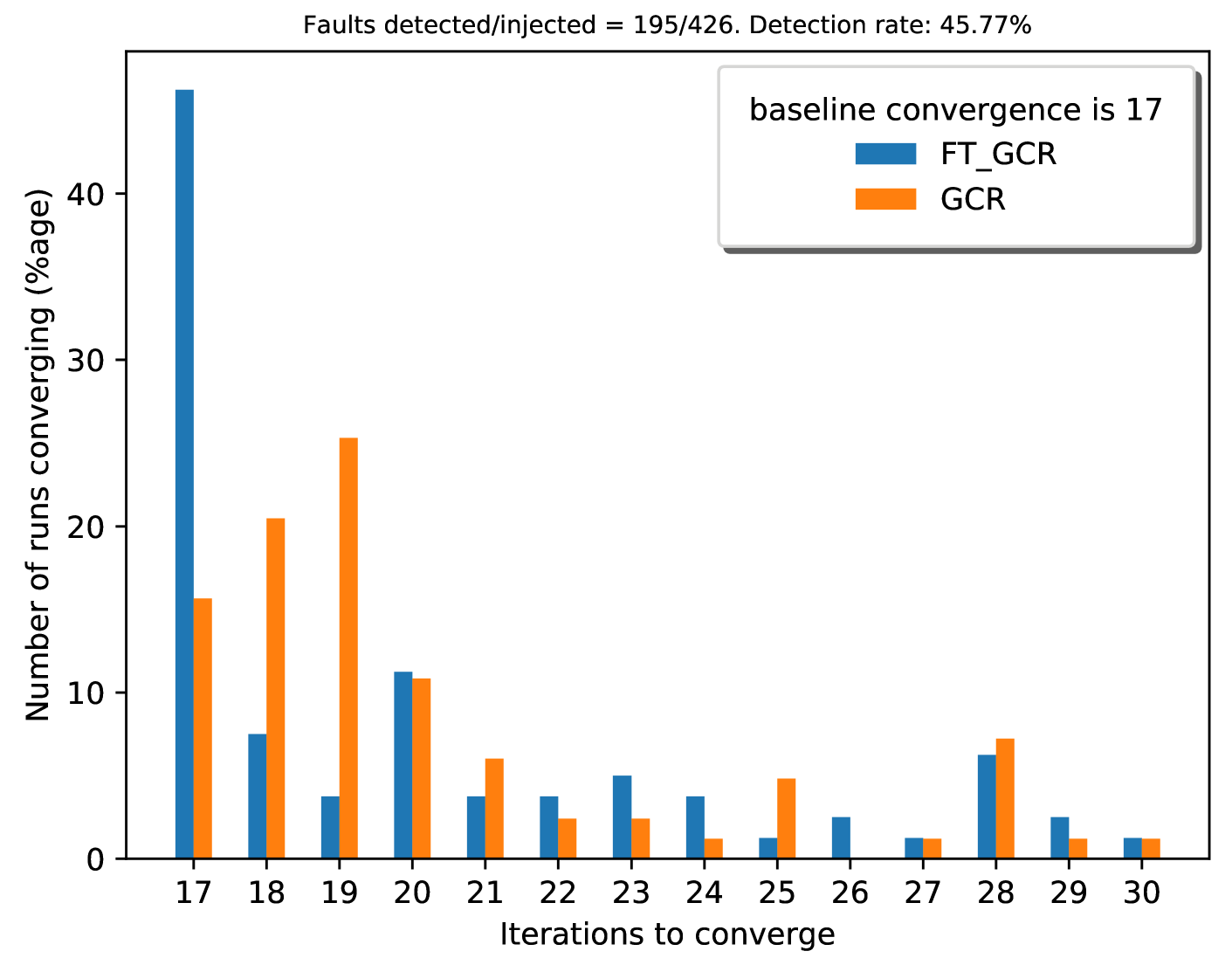}
}\\
\subfloat{
\includegraphics[width=0.46 \textwidth]{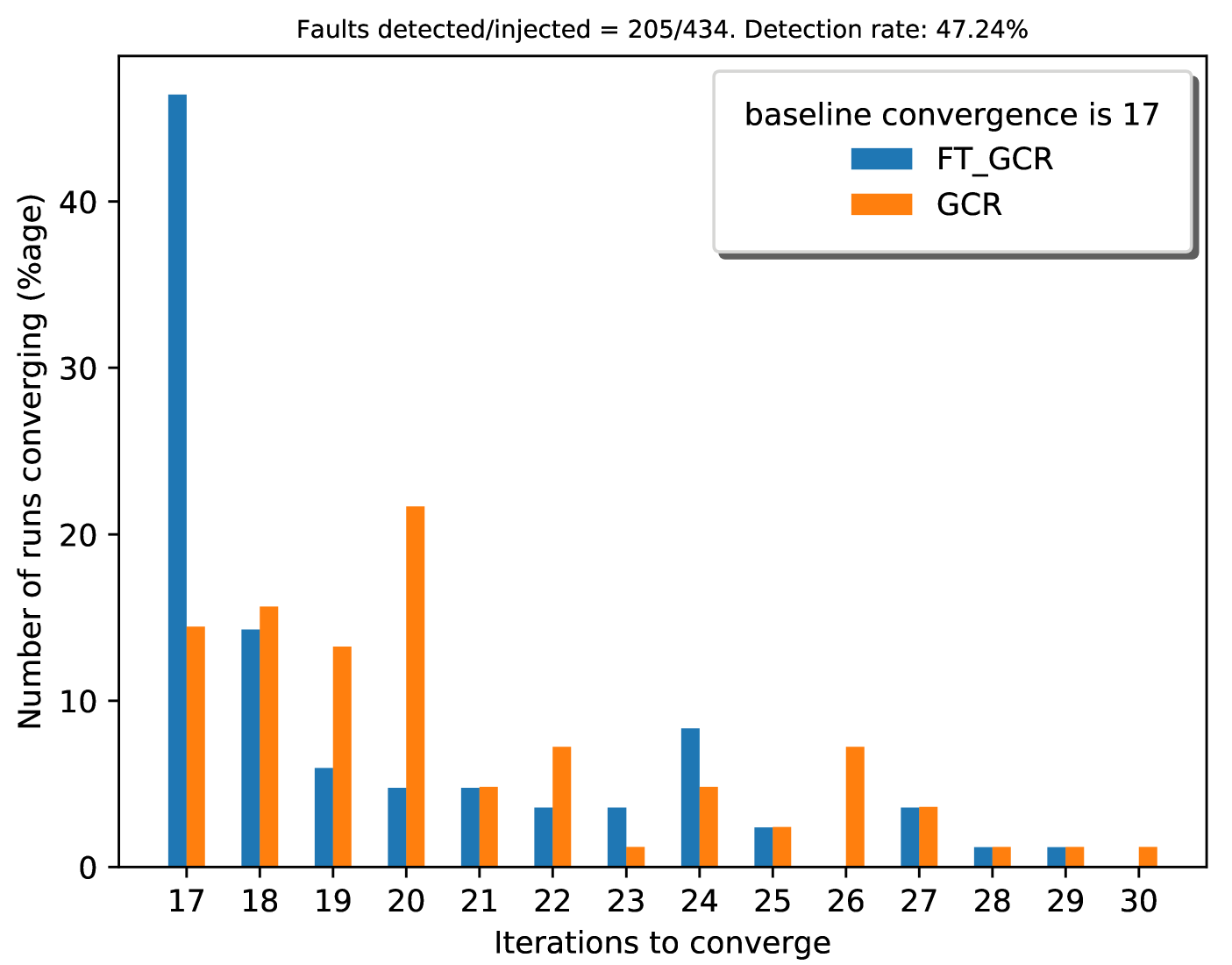}
}\subfloat{
\includegraphics[width=0.46 \textwidth]{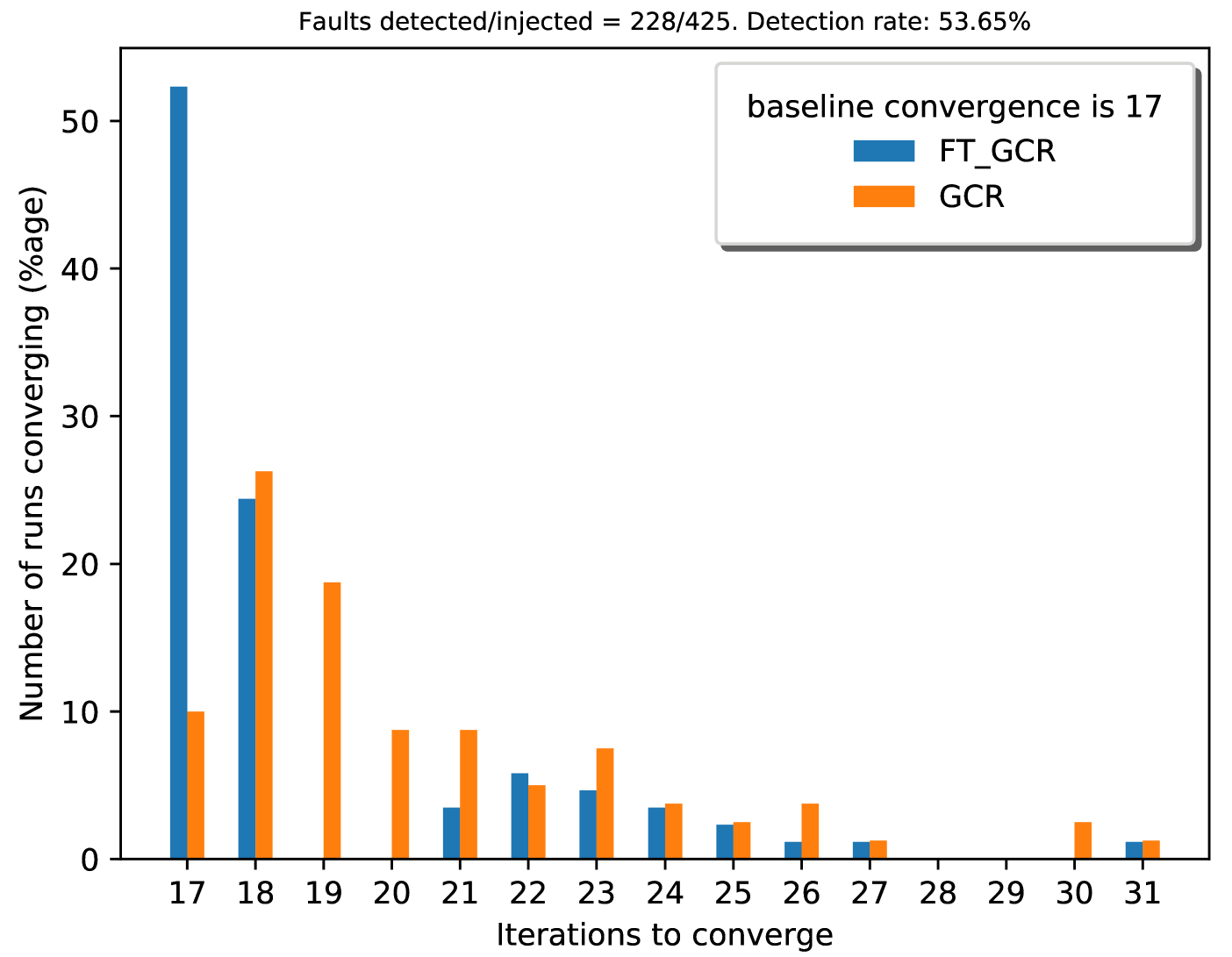}
}
   
\caption{Numerical behaviour of GCR and FT-GCR when faults are introduced, O80 grid, 72 MPI tasks, size 5 Krylov space, 2\% fault injection probability. Each histogram displays the percentage of runs that converged at a given iteration. Top left: 0.0004\% data loss, top right: 0.04\% data loss, middle left: 0.2\% data loss, middle right: 1\% data loss, bottom left: 5\% data loss, bottom right: 20\% data loss.}
\label{fig:GCRconvergence_O80}
\end{figure}

Figure \ref{fig:GCRconvergence_O80} demonstrates results for the runs with O80 grid, $2\%$ fault injection probability, see also Table \ref{tbftgcr_O40_O640} in the Appendix \ref{appendix} for results with the range of grids O40 to O640. Each output dataset features results from around $100$ protected runs and $100$ unprotected runs, both batches using the same parameters, but with fully randomised faults (the seeds of both fault events and the data values corrupted are unique for each run\footnote{Note that actual data presented will be less due to a $\sim 20\%$ share of non-faulty runs, which are disregarded}). Baseline convergence of an un-corrupted simulation is also established by running without fault injection. 

In evaluating the results, it should first be noted that if all faults are detected by FT-GCR, convergence will match the baseline. This does not mean that there is no computational cost to the detection and correction of such faults -- any calculations performed between the backup of a known good configuration and the detection of a fault are wasted. These calculations amount to, at most, one full GCR iteration (denoted by $n$ in Algorithm \ref{ftgcrk}) per fault. In addition, undetected faults may not always have an appreciable effect on the convergence. Upon examination of individual runs, early faults -- those injected in the first few iterations -- often go undetected, the reason being the sharp reduction in solution error in the first iterations of the elliptic solver. Therefore, errors caused by early faults are usually hidden by the overall reduction in the residual norm, yet those same faults often appear to have more significant implications to stability than faults injected towards the end of the simulation, where the error is much smoother. Late faults also tend to be detected more often, and the solver seems to fare better when the general solution error is already both small and smooth. 

For small fault events (e.g., at $0.0004\%$ data loss), faults are both difficult to detect and have a low impact on convergence, indicating that at least for this test case GCR is highly resilient to soft faults. In fact, even for very implausibly large fault events (e.g.\ $5\%$ or $20\%$) convergence is only delayed by at most $\sim\!20\%$. These large corruption events are, in effect, altering the values of every data entry in the entire array on a given process. That this type of event does not break the overall solver -- even occasionally -- is quite remarkable. The results also indicate that detection of data corruption becomes increasingly more difficult as the number of grid nodes increases (even when corruption events are scaled by this number).

\newpage 

\begin{figure}
    \centering
    \includegraphics[width=\columnwidth]{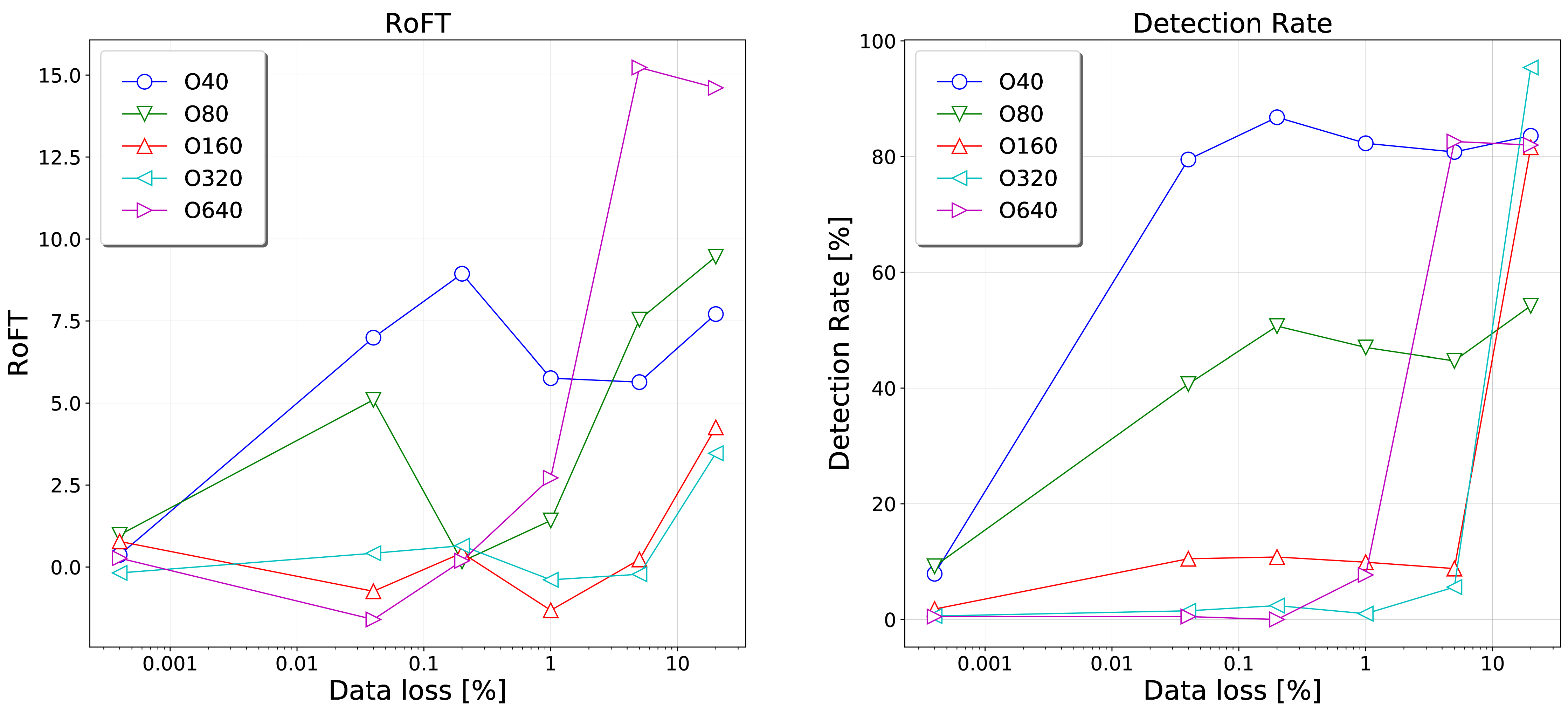}
    \caption{Return on Fault Tolerance score (RoFT, left) and FT-GCR Detection rate (right) as a function of data loss for the tests in Table \ref{tbftgcr_O40_O640}.}
    \label{fig:Roft_Detec_rate_vs_frequency}
\end{figure}

Figure \ref{fig:Roft_Detec_rate_vs_frequency} displays the relationship between the size of a fault event (labelled as data loss), the improvement in average convergence (left), and the detection rate of faults (right).
The Return on Fault Tolerance (RoFT) score for a grid and percentage of data lost is computed as:
\begin{equation}
\textrm{RoFT} = \dfrac{\textrm{GCR} - \textrm{FT-GCR}}{\textrm{baseline}}    
\end{equation}
where `GCR' and `FT-GCR' denote the average number of iterations to convergence for unprotected runs  and protected runs, respectively, normalised by the baseline convergence. The RoFT value intends to represent the savings in convergence delay when the model is protected using the FT-GCR algorithm compared to the basic GCR algorithm. As noted above, the benefits of FT-GCR tend to strongly depend on the overall size of individual fault events. Larger faults are better detected, and detection of smaller faults is improved on smaller grids.

\begin{figure}
    \centering
    \includegraphics[width=.6\columnwidth]{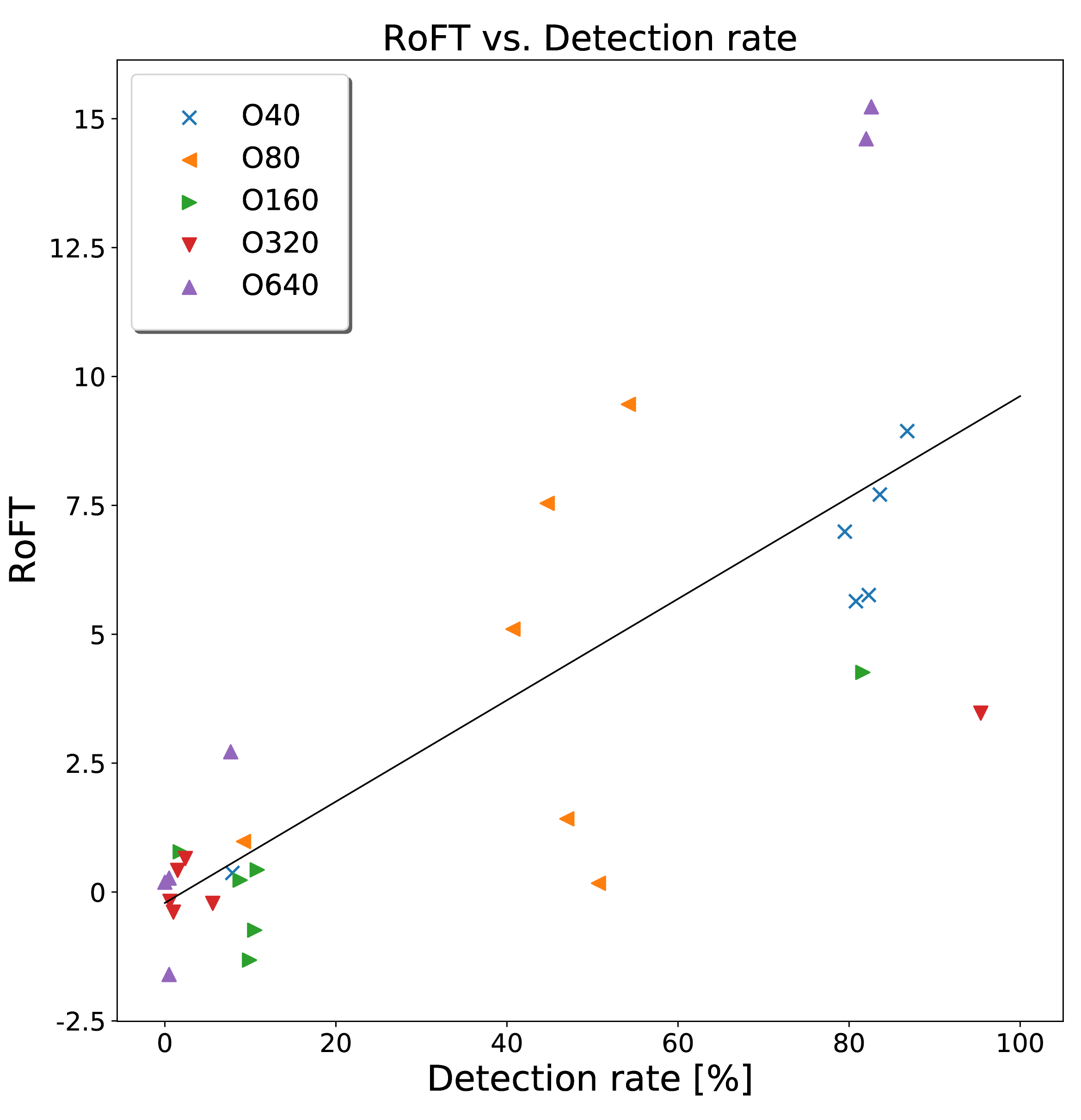}
    \caption{Scatter plot of Return on Fault Tolerance (RoFT) and FT-GCR Detection rate values and linear fit ($R^2=0.608$) for the tests in Table \ref{tbftgcr_O40_O640}.}
    \label{fig:Roft_Detec_rate_scatter}
\end{figure}

Figure \ref{fig:Roft_Detec_rate_scatter} demonstrates the relationship between detection rates and improved convergence. This relationship is not directly linear, but the overall trend is that better detection has a positive effect on convergence. The data towards the bottom left of the plot is less significant, since at low detection rates the unknown effect of a fault on convergence is difficult to quantify. The expectation is that more impactful faults will be more often detected, while the less impactful faults are more likely to be missed. Limited testing indicates this is not a reliable assumption, for example faults early in the run are detected with lower data loss values.

In addition, increasing the number of faults per run using a 5\% fault injection probability both increases the convergence delay of unprotected runs and increases the detection rate of faults (Table \ref{tbftgcr_O160_O320_5pc_inj_prob} in Appendix \ref{appendix}). Figure \ref{fig:Roft_Detec_rate_vs_frequency_injsens} displays the corresponding RoFT and detection rates vs.\ the relative size of the fault event for the simulations with the O$160$ and O$320$ grids, contrasting the runs at $2\%$ and $5\%$ fault injection probability. In both cases, FT-GCR provides an improvement in performance as the number of fault events increases, with a more pronounced effect seen for very large fault events ($5\%$ and $20\%$ data loss values). 

It is clear that when most faults are detected, FT-GCR provides a noticeable improvement in convergence characteristics in comparison to unprotected GCR simulations. Whenever this is not the case, the impact of fault tolerance is more limited. 

\begin{figure}[H]
    \centering
    \includegraphics[width=.9\columnwidth]{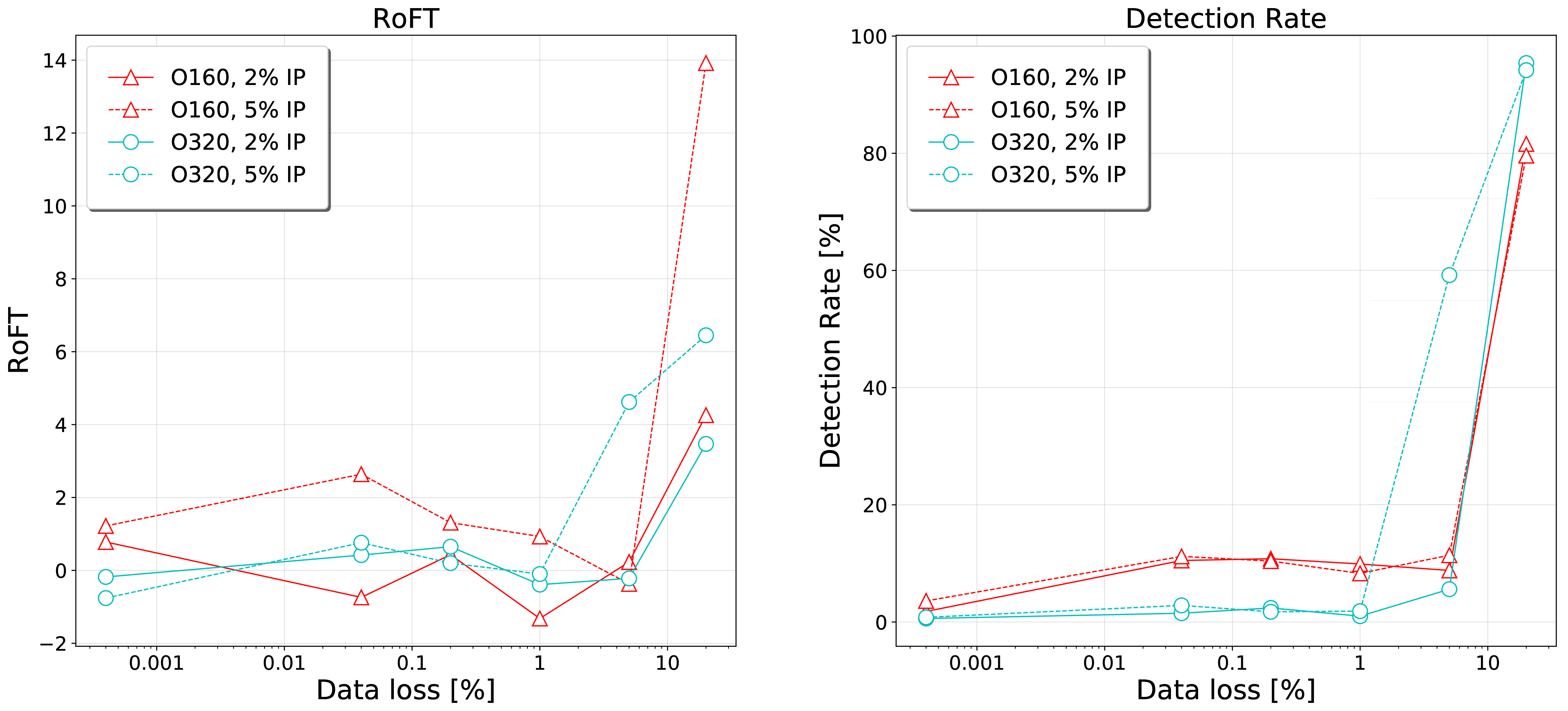}
    \caption{Return on Fault Tolerance score (RoFT, left) and FT-GCR Detection rate (right) as a function of data loss, O160 and O320 grids, 2\% vs.\ 5\% fault injection probability sensitivity tests in Table \ref{tbftgcr_O160_O320_5pc_inj_prob}.}
    \label{fig:Roft_Detec_rate_vs_frequency_injsens}
\end{figure}

\section{Discussion and Conclusion}

We have introduced FT-GCR, a fault-tolerant version of the preconditioned Generalized Conjugate Residual elliptic solver. The fault detection mechanism of FT-GCR is based on a simple check on the residual norm decrease at each iteration. Once a fault is detected, recovery from the fault is based on a restart of the computation within the current Krylov subspace. Although the detection process uses a quantity whose evaluation requires global communications, it is also a quantity that is already computed. The backup of a known-good configuration increases storage requirements of the GCR implementation by a small amount -- $4$ out of $14$ full mesh arrays are backed up, corresponding to $\approx30\%$ memory overhead for the GCR routine. Implementation within an existing GCR solver is therefore very simple, requiring only the addition of a backup loop, a fault detection subroutine, and a restart loop. The main workings of the algorithm remain unchanged. 

The modified solver has been tested against an unmodified version on the solution of steady potential flow over a hill on the planet using the Octahedral Gaussian grid at horizontal resolutions ranging from $278\,\textrm{km}$ to $15.5\,\textrm{km}$ - therefore up to about a half of the current operational resolution of state-of-the-art medium-range weather forecast models - and considering a variety of choices for fault parameters such as percentage of data corrupted and chance of fault injection.

Parallel simulations with up to 3240 MPI processes have validated the capabilities of the method in protecting the solver's performance when faults in the form of bit flips are introduced. The rate of detection of faults by FT-GCR has been shown to grow with the percentage of data corrupted, and when the fault-tolerant version is used, runs converge in fewer iterations than with standard GCR in the vast majority of cases. The return on fault tolerance is larger for higher detection rates, and both parameters increase with increasing resolution, an effect best seen with the finest grids and largest proportions of data corrupted by the faults.

Besides the overall picture of a positive effect of the modifications introduced for algorithmic resilience, the experiments have also highlighted a substantial degree of inherent robustness in the standard version of the GCR method. Even with massive fault injections, most GCR runs displayed iteration counts very close both to the faultless baseline run and to the FT-GCR run. Indeed no situations have been found where faults caused permanent stagnation of runs. 
It should also be noted that the variation of FT-GCR and GCR plots, as well as the negative values for the return parameter, at low detection rates are a result of the variability of runs. It would likely require thousands of additional data points to properly display the small improvement brought by low detection rates.

The tests presented here have been performed across a range of resolutions for demonstration purposes. 
The convergence acceleration and detection rate associated with FT-GCR were most pronounced at very large values of data loss, while the skill scores of the method for smaller, more realistic values of data loss seem less conclusive in the present version. Therefore, as currently implemented the method would be of limited use in modelling scenarios where it is important to achieve large detection rates even with very small corruption events. At the same time, the low computational cost and ease of implementation within a standard Krylov method make FT-GCR a useful tool in numerical frameworks where such solvers can take up a significant portion of wallclock time. In the context of parallel runs using FT-GCR as linear solver, for example within semi-implicit time discretizations in computational models of fluid flow, the fault tolerance against bit flips provided by a method such as the one presented here would ideally be complemented with recovery from node failures and systems resilience in an overall strategy as discussed in \cite{benacchio:2021, giraud:2020}.

\subsection*{Acknowledgments}
The authors were supported by the ESCAPE-2 project, European Union’s Horizon 2020 research and innovation programme (grant agreement No 800897).

\appendix
\section{Complete results for the runs}
\label{appendix}

\small 
Tables \ref{tbftgcr_O40_O640} and \ref{tbftgcr_O160_O320_5pc_inj_prob} contain the complete results of the tests performed with the FT-GCR fault tolerant solver to produce  Figures \ref{fig:Roft_Detec_rate_vs_frequency}, 
\ref{fig:Roft_Detec_rate_scatter} and \ref{fig:Roft_Detec_rate_vs_frequency_injsens} in the main text.

\begin{table}[H]\setlength{\tabcolsep}{3pt}
\scriptsize
\caption{Average number over all runs of faults per run and faults detected, detection rate for FT-GCR, average iterations for FT-GCR, average iterations for unprotected GCR, and Return on Fault Tolerance for each run with the grids and data loss values of Table \ref{tab:resol}, 2\% fault injection probability. In the baseline configuration, convergence is reached in  $\{19,17,17,21,21\}$ iterations for the  $\{\textrm{O}40,\textrm{O}80,\textrm{O}160,\textrm{O}320,\textrm{O}640\}$ grids.}
\vspace{\baselineskip}
\centering
\begin{tabular}{lrrrrrrr}
\toprule
\multirow{2}{*}{Grid} & \multirow{2}{*}{$\overset{\textrm{\small $\%$ data}}{\quad\,\text{\small loss}}$} & \multirow{2}{*}{$\overset{\quad\textrm{\small Faults}}{\;\,\underset{\text{\small per run}}{\large\textcolor{white}{.}}}$} & \multirow{2}{*}{$\overset{\quad\textrm{\small Faults}}{\text{\small detected}}$} & \multirow{2}{*}{$\overset{\quad\textrm{\small Detection}}{\text{\small \quad Rate [$\%$]}}$} &   \multirow{2}{*}{$\overset{\textrm{\small Convergence}}{\quad\;\text{\small FT-GCR}}$} & \multirow{2}{*}{$\overset{\textrm{\small Convergence}}{\qquad\quad\text{\small GCR}}$} & \multirow{2}{*}{$\overset{\textrm{\small RoFT}}{\quad\;\text{\small [$\%$]}}$} \\
& & & & & & \\
\midrule
\multirow{6}{*}{O40}  & $0.0004$ & 2.56 & 0.17 & 7.9 & 19.1 & 19.16 & 0.37 \\
                      & $0.04$   & 2.80 & 2.20 & 79.5 & 19.15 & 20.48 & 6.99 \\
                      & $0.2$    & 2.94 & 2.43 & 86.8 & 19.17 & 20.87 & 8.94 \\
                      & $1.0$    & 3.02 & 2.49 & 82.3 & 19.28 & 20.37 & 5.76 \\
                      & $5.0$    & 2.58 & 2.11 & 80.8 & 19.23 & 20.31 & 5.64 \\
                      & $20.0$   & 2.89 & 2.45 & 83.6 & 19.15 & 20.62 & 7.71 \\
\midrule
\multirow{6}{*}{O80}  & $0.0004$ & 2.36 & 0.3 & 9.2 & 18.83 & 19    & 0.98 \\
                      & $0.04$   & 2.4  & 1.02  & 40.7 & 19.46 & 20.33 & 5.1 \\
                      & $0.2$    & 2.79 & 1.33 & 50.7 & 20.32 & 20.35 & 0.17 \\
                      & $1.0$    & 2.61 & 1.2 & 47 & 19.99 & 20.23 & 1.42 \\
                      & $5.0$    & 2.6 & 1.23 & 44.7 & 19.39 & 20.67 & 7.54 \\
                      & $20.0$   & 2.56 & 1.37 & 54.2 & 18.77 & 20.38 & 9.46 \\
\midrule
\multirow{6}{*}{O160}  & $0.0004$ & 2.15 & 0.06 & 1.8 & 17.80 & 17.94 & 0.78  \\
                       & $0.04$   & 2.27 & 0.18 & 10.5 & 17.93 & 17.8  & -0.74 \\
                       & $0.2$    & 2.04 & 0.21 & 10.8 & 17.55 & 17.62 & 0.43  \\
                       & $1.0$    & 2.23 & 0.2 & 9.9 & 17.59 & 17.36 & -1.32 \\
                       & $5.0$    & 2.15 & 0.19 & 8.8 & 17.74 & 17.78 & 0.23  \\
                       & $20.0$   & 2.19 & 1.72 & 81.6 & 17.4    & 18.13 & 4.26  \\
\midrule
\multirow{6}{*}{O320}  & $0.0004$ & 2.23 & 0.01 & 0.6  & 21.08 & 21.04 & -0.18      \\
                       & $0.04$   & 2.28 & 0.04 & 1.5 & 21 & 21.09 & 0.42     \\
                       & $0.2$    & 2.48 & 0.05 & 2.4  & 21.1 & 21.23 & 0.65     \\
                       & $1.0$    & 2.37 & 0.02 & 1  & 21.09 & 21.01 & -0.39     \\
                       & $5.0$    & 2.36 & 0.1 & 5.6  & 21.05 & 21.01 & -0.22     \\
                       & $20.0$   & 2.58 & 2.5 & 95.4  & 20.99 & 21.72 & 3.47
\\
\midrule
\multirow{6}{*}{O640}  & $0.0004$ & 2.33 & 0.02 & 0.5 & 21.13 & 21.18 & 0.27     \\
                       & $0.04$   & 2.4 & 0.01 & 0.5  & 21.44 & 21.10 & -1.6     \\
                       & $0.2$    & 2.26 & 0.02 & 0  & 21.37 & 21.41 & 0.2     \\
                       & $1.0$    & 2.81 & 0.28 & 7.7  & 23.86 & 24.43 & 2.72     \\
                       & $5.0$    & 2.70 & 2.26 & 82.6  & 22.1 & 25.29 & 15.23      \\
                       & $20.0$   & 2.42 & 2.1 & 82 & 22.39 & 25.46 & 14.61
\\
\bottomrule
\end{tabular}
\vspace{\baselineskip}
\label{tbftgcr_O40_O640}
\end{table}

\begin{table}[H]\setlength{\tabcolsep}{3pt}
\scriptsize
\caption{As in Table \ref{tbftgcr_O40_O640} but for the O160, O320 grids only and with 5\% fault injection probability. In the baseline configuration, convergence is reached in  $\{17,21\}$ iterations for the  $\{\textrm{O}160,\textrm{O}320\}$ grids.}
\vspace{\baselineskip}
\centering
\begin{tabular}{lrrrrrrr}
\toprule
\multirow{2}{*}{Grid} & \multirow{2}{*}{$\overset{\textrm{\small $\%$ data}}{\quad\,\text{\small loss}}$} & \multirow{2}{*}{$\overset{\quad\textrm{\small Faults}}{\;\,\underset{\text{\small per run}}{\large\textcolor{white}{.}}}$} & \multirow{2}{*}{$\overset{\quad\textrm{\small Faults}}{\text{\small detected}}$} & \multirow{2}{*}{$\overset{\quad\textrm{\small Detection}}{\text{\small \quad Rate [$\%$]}}$} &   \multirow{2}{*}{$\overset{\textrm{\small Convergence}}{\quad\;\text{\small FT-GCR}}$} & \multirow{2}{*}{$\overset{\textrm{\small Convergence}}{\qquad\quad\text{\small GCR}}$} & \multirow{2}{*}{$\overset{\textrm{\small RoFT}}{\quad\;\text{\small [$\%$]}}$} \\
& & & & & & \\
\midrule
\multirow{6}{*}{O160}  & $0.0004$ & 4.73 & 0.14 & 3.6 & 17.97 & 18.18 & 1.22  \\
                       & $0.04$   & 4.61 & 0.51 & 11.2 & 18.15 & 18.6  & 2.64 \\
                       & $0.2$    & 4.71 & 0.47 & 10.4 & 18.79 & 19.01 & 1.31  \\
                       & $1.0$    & 4.71 & 0.52 & 8.3 & 18.25 & 18.41 & 0.93 \\
                       & $5.0$    & 4.6 & 0.57 & 11.4 & 18.41 & 18.35 & -0.37  \\
                       & $20.0$   & 5.36 & 4.14 & 79.6 & 17.91 & 20.27 & 13.92  \\
\midrule
\multirow{6}{*}{O320}  & $0.0004$ & 4.86 & 0.04 & 0.8  & 21.23 & 21.07 & -0.76      \\
                       & $0.04$   & 4.73 & 0.1 & 2.84 & 21.13 & 21.29 & 0.76     \\
                       & $0.2$    & 5.35 & 0.09 & 1.74  & 21.12 & 21.16 & 0.2     \\
                       & $1.0$    & 5.13 & 0.09 & 1.86  & 21.2 & 21.18 & -0.1     \\
                       & $5.0$    & 5.61 & 3.16 & 59.21  & 21.42 & 22.39 & 4.62     \\
                       & $20.0$   & 5.74 & 5.49 & 94.18  & 21.02 & 22.37 & 6.45
\\
\bottomrule
\end{tabular}
\vspace{\baselineskip}
\label{tbftgcr_O160_O320_5pc_inj_prob}
\end{table}

\section*{Note}

This work has not yet been peer-reviewed and is provided by the contributing authors as a means to ensure timely dissemination of scholarly and technical work on a noncommercial basis. Copyright and all rights therein are maintained by the authors or by other copyright owners. It is understood that all persons copying this information will adhere to the terms and constraints invoked by each author's copyright. This work may not be reposted without explicit permission of the copyright owner.


\footnotesize

\bibliographystyle{plain}
\bibliography{Bibliography}


\end{document}